\newcount\secno
\newcount\prmno
\newif\ifnotfound
\newif\iffound

\def\section#1{\vskip1truecm
               \global\def\currenvir{section}
               \global\advance\secno by1\global\prmno=0
               {\bf \number\secno. {#1}}
               \vglue2pt}

\def\subsection{\global\def\currenvir{subsection}
                \global\advance\prmno by1
               \medskip
 \ind{ (\number\secno.\number\prmno) }}
\def\subsec{\global\def\currenvir{subsection}
                \global\advance\prmno by1
                { (\number\secno.\number\prmno)\ }}

\def\proclaim#1{\global\advance\prmno by 1
                {\bf #1 \the\secno.\the\prmno$.-$ }}

\long\def\th#1 \enonce#2\endth{%
   \vglue6pt\proclaim{#1}{\it #2}\global\def\currenvir{th}\smallskip}

\def\bib#1{\rm #1}
\long\def\thr#1\bib#2\enonce#3\endth{%
\medbreak{\global\advance\prmno by 1\bf#1\the\secno.\the\prmno\ 
\bib{#2}$\!.-$ } {\it
#3}\global\def\currenvir{th}\smallskip}
\def\rem#1{\global\advance\prmno by 1
{\it #1} \the\secno.\the\prmno$.-$ }


\magnification 1200
\pretolerance=500 \tolerance=1000  \brokenpenalty=5000
\mathcode`A="7041 \mathcode`B="7042 \mathcode`C="7043
\mathcode`D="7044 \mathcode`E="7045 \mathcode`F="7046
\mathcode`G="7047 \mathcode`H="7048 \mathcode`I="7049
\mathcode`J="704A \mathcode`K="704B \mathcode`L="704C
\mathcode`M="704D \mathcode`N="704E \mathcode`O="704F
\mathcode`P="7050 \mathcode`Q="7051 \mathcode`R="7052
\mathcode`S="7053 \mathcode`T="7054 \mathcode`U="7055
\mathcode`V="7056 \mathcode`W="7057 \mathcode`X="7058
\mathcode`Y="7059 \mathcode`Z="705A
\def\spacedmath#1{\def\packedmath##1${\bgroup\mathsurround =0pt##1\egroup$}
\mathsurround#1
\everymath={\packedmath}\everydisplay={\mathsurround=0pt}}
\def\nospacedmath{\mathsurround=0pt
\everymath={}\everydisplay={} } \spacedmath{2pt}
\def\qfl#1{\buildrel {#1}\over {\longrightarrow}}
\def\phfl#1#2{\normalbaselines{\baselineskip=0pt
\lineskip=10truept\lineskiplimit=1truept}\nospacedmath\smash {\mathop{\hbox to
8truemm{\rightarrowfill}}
\limits^{\scriptstyle#1}_{\scriptstyle#2}}}
\def\hfl#1#2{\normalbaselines{\baselineskip=0truept
\lineskip=10truept\lineskiplimit=1truept}\nospacedmath\smash{\mathop{\hbox to
12truemm{\rightarrowfill}}\limits^{\scriptstyle#1}_{\scriptstyle#2}}}

\def\iso{\vbox{\hbox to .8cm{\hfill{$\scriptstyle\sim$}\hfill}
\nointerlineskip\hbox to .8cm{{\hfill$\longrightarrow $\hfill}} }}

\def\sdir_#1{\mathrel{\mathop{\kern0pt\oplus}\limits_{#1}}}
\def\pprod_#1^#2{\raise
2pt \hbox{$\mathrel{\scriptstyle\mathop{\kern0pt\prod}\limits_{#1}^{#2}}$}}
\def\bot_#1{\mathrel{\mathop{\kern0pt\bigotimes}\limits_{#1}}}
 \newfam\bboardfam
\font\eightrm=cmr8         \font\eighti=cmmi8
\font\eightsy=cmsy8        \font\eightbf=cmbx8
\font\eighttt=cmtt8        \font\eightit=cmti8
\font\eightsl=cmsl8        \font\sixrm=cmr6
\font\sixi=cmmi6           \font\sixsy=cmsy6
\font\sixbf=cmbx6
\font\eightbboard=msbm10 at 8pt\font\sevenbboard=msbm7
\catcode`\@=11
\def\eightpoint{%
  \textfont0=\eightrm \scriptfont0=\sixrm \scriptscriptfont0=\fiverm
  \def\rm{\fam\z@\eightrm}%
  \textfont1=\eighti  \scriptfont1=\sixi  \scriptscriptfont1=\fivei
  \def\oldstyle{\fam\@ne\eighti}\let\old=\oldstyle
  \textfont2=\eightsy \scriptfont2=\sixsy \scriptscriptfont2=\fivesy
  \textfont\itfam=\eightit
  \def\it{\fam\itfam\eightit}%
  \textfont\slfam=\eightsl
  \def\sl{\fam\slfam\eightsl}%
  \textfont\bffam=\eightbf \scriptfont\bffam=\sixbf
  \scriptscriptfont\bffam=\fivebf
  \def\bf{\fam\bffam\eightbf}%
  \textfont\ttfam=\eighttt
  \def\tt{\fam\ttfam\eighttt}%
  \textfont\bboardfam=\eightbboard \scriptfont\bboardfam=\sevenbboard
  \scriptscriptfont\bboardfam=\sevenbboard
  \def\bb{\fam\bboardfam}%
  
  \abovedisplayskip=9pt plus 3pt minus 9pt
  \belowdisplayskip=\abovedisplayskip
  \abovedisplayshortskip=0pt plus 3pt
  \belowdisplayshortskip=3pt plus 3pt 
  \smallskipamount=2pt plus 1pt minus 1pt
  \medskipamount=4pt plus 2pt minus 1pt
  \bigskipamount=9pt plus 3pt minus 3pt
  \normalbaselineskip=9pt
  \setbox\strutbox=\hbox{\vrule height7pt depth2pt width0pt}%
  \normalbaselines\rm}\catcode`\@=12

\newcount\noteno
\noteno=0
\def\pc#1{\tenrm#1\sevenrm}
\def\tx{\kern-1.5pt -}
\def\cqfd{\kern 2truemm\unskip\penalty 500\vrule height 4pt depth 0pt width
4pt\medbreak} 
\def\virg{\raise
.4ex\hbox{,}}
\def\ind{\par\hskip 0.8truecm\relax}

\def\rond{\kern 1pt{\scriptstyle\circ}\kern 1pt}

\def\Hom{\mathop{\rm Hom}\nolimits}
\def\Ext{\mathop{\rm Ext}\nolimits}

\def\im{\mathop{\rm Im}\nolimits}
\def\Ker{\mathop{\rm Ker}\nolimits}

\def\Pic{\mathop{\rm Pic}\nolimits}
\def\Div{\mathop{\rm Div}\nolimits}

\def\div{\mathop{\rm div\,}\nolimits}
\def\Br{\mathop{\rm Br}\nolimits}

\def\bb{\Bbb}
\def\Z{\bb Z}
\def\Q{\bb Q}
\def\R{\bb R}
\def\C{\bb C}
\def\P{\bb P}

\def\gm{{\bb G}_m}
\def\Nm{\mathop{\rm Nm}\nolimits}

\input amssym.def
\input amssym
\input xy
\xyoption{all}

\font\san=cmssdc10

\def\sym{\hbox{\san \char83}}

\frenchspacing
\catcode`\@=11\def\eqalign#1{\null\,\vcenter{\openup\jot\m@th\ialign{
\strut\hfil$\displaystyle{##}$&$\displaystyle{{}##}$
&&\quad\strut$\displaystyle{##}$&$\displaystyle{{}##}$
\crcr#1\crcr}}\,}
\catcode`\@=12

\vsize = 23.5truecm
\hsize = 16truecm
\hoffset =0.5truecm
\voffset = -0.5truecm

\parindent=0cm
\baselineskip15pt
\overfullrule=0pt

\centerline{\bf  On the Brauer group of  Enriques surfaces}\smallskip
 \centerline{Arnaud {\pc BEAUVILLE}} 
\vskip.7cm

\section{Introduction}

\ind Let $S$ be a complex Enriques surface, and $\pi :X\rightarrow S$ its 2-to-1 cover by a K3 surface. Poincar\'e duality provides an isomorphism $H^3(S,\Z)\cong H_1(S,\Z)=\Z/2$, so that there is a unique nontrivial element $b_S$ in the Brauer group  $\Br(S)$. What is the pull-back of this element in $\Br(X)$? Is it nonzero?
\ind The answer to the first question is easy in terms of the canonical isomorphism $\Br(X)\iso \Hom(T_X,\Q/\Z)$ (see \S 2): $\pi ^*b_S$ corresponds to the linear form $\tau \mapsto (\beta \cdot \pi _*\tau )$, where $\beta $ is any element of $H^2(S,\Z/2)$ which does not come from $H^2(S,\Z)$. 
The second question turns out to be more subtle: the answer depends on the surface. We will characterize the  surfaces $S$ for which $\pi ^*b_S=0$ (Corollary 5.7), and show that they form a countable union of hypersurfaces in the moduli space of
Enriques surfaces (Corollary 6.5). 
\ind Part of our results hold over any algebraically closed field, and also in a more general set-up (see Proposition 4.1 below); for the last part, however, we need in a crucial way Horikawa's description of the moduli space by transcendental methods. 

 {\bigskip {\eightrm\baselineskip=12pt
\leftskip1cm\rightskip1cm\hskip0.8truecm  The question considered here is mentioned in [H-S], Cor. 2.8. I am  indebted to
\hbox{J.-L. Colliot-Th\'el\`ene}  for explaining it to me, and for very
useful discussions and comments. I am  grateful to J.~Lannes
for providing the topological proof of Lemma 5.4.\par}}\smallskip  

\section{The Brauer group of a surface}
\ind Let $S$ be a smooth projective variety over a field; we define  the Brauer group $\Br(S)$ as the \'etale cohomology group $H^2_{\rm{\acute e t}}(S,\gm)$. This definition coincides with that of Grothendieck [G] by a result of Gabber, which we will not need here.

\ind In this section we assume  that $S$ is a complex  surface; we  recall the  description of  $\Br(S)$ in that case  -- this is classical but not so easy to find in the literature. The Kummer exact sequence
$$0\rightarrow \Z/n\longrightarrow \gm\longrightarrow \gm \rightarrow 0$$gives rise to an 
exact sequence
$$0\rightarrow \Pic(S)\otimes \Z/n \longrightarrow H^2(S,\Z/n)\qfl{p}   \Br(S)[n]\rightarrow 0\eqno(\number\secno.a)$$
(we denote by $M[n]$ the kernel of the multiplication by $n$ in a  $\Z$\tx module $M$).

\ind On the other hand,  the cohomology  exact sequence associated to \hfill\break $0\rightarrow \Z\rightarrow \Z\rightarrow \Z/n\rightarrow 0$ gives:
$$0\rightarrow H^2(S,\Z)\otimes \Z/n \longrightarrow H^2(S,\Z/n)\longrightarrow H^3(S,\Z)[n]\rightarrow 0 \eqno(\number\secno.b)$$
Comparing $(\number\secno.a)$ and $(\number\secno.b)$  we get an  exact sequence
$$0\rightarrow  \Pic(S)\otimes \Z/n \longrightarrow H^2(S,\Z)\otimes \Z/n\longrightarrow \Br(S)[n]\longrightarrow  H^3(S,\Z)[n]\rightarrow 0\ .\eqno(\number\secno.c)$$

Let $H^2(S,\Z)_{\rm tf}$ be the quotient of $H^2(S,\Z)$ by its torsion subgroup; the cup-product induces a perfect pairing on $H^2(S,\Z)_{\rm tf}$. We denote by 
 $T_S\subset H^2(S,\Z)_{\rm tf}$  the {\it transcendental lattice}, that is, the orthogonal of the image of $\Pic(S)$. We have an exact sequence
$$ \Pic(S)\qfl{c_1} H^2(S,\Z)\qfl{u} T_S^*\rightarrow 0$$
where $u$ associates to $\alpha \in H^2(S,\Z)$ the cup-product with $\alpha $.
Taking tensor product with $\Z/n$ and comparing with $(\number\secno.c)$, we get 
 an exact sequence
$$0\rightarrow \Hom(T_S,\Z/n)\longrightarrow \Br(S)[n]\longrightarrow  H^3(S,\Z)[n]\rightarrow 0\ ; \eqno(\number\secno.d)$$or, passing to the direct limit over $n$,
$$0\rightarrow \Hom(T_S,\Q/\Z)\longrightarrow \Br(S)\longrightarrow  {\rm Tors}\,H^3(S,\Z)\rightarrow 0\ . \eqno(\number\secno.e)$$

\section{Topology of Enriques surfaces}
\subsection  Let $S$ be an Enriques surface  (over $\C$). 
We first recall some elementary facts on the topology of $S$. A general reference is [BHPV], ch. VIII.

\ind  The torsion subgroup of  $H^2(S,\Z)$ is isomorphic to $\Z/2$; its
 nonzero element is the canonical class $K_S$. 
  Let $k_S$ denote the image of $K_S$ in $H^2(S,\Z/2)$. The universal coefficient theorem together with Poincar\'e duality gives an exact sequence
$$0\rightarrow \Z/2 \qfl{k_S} H^2(S,\Z/2)\qfl{v_S} \Hom(H^2(S,\Z),\Z/2)\rightarrow 0 \eqno(\number\secno.a)$$where $v_S$ is deduced from the cup-product. 
\subsection  The linear form $\alpha \mapsto (k_S\cdot \alpha )$ on $H^2(S,\Z/2)$ vanishes on the image of $H^2(S,\Z)$, hence coincides with the map 
  $H^2(S,\Z/2)\rightarrow H^3(S,\Z)=\Z/2$ from the
exact sequence $(2.b)$. Note  that $k_S$ is  the second Stiefel-Whitney class $w_2(S)$; in particular, we have $( k_S \cdot\alpha )=\alpha ^2$ for all $\alpha \in H^2(S,\Z/2)$ (Wu formula, see [M-S]).
\subsection The map $c_1:\Pic(S)\rightarrow H^2(S,\Z)$ is an isomorphism, hence $(2.e)$ provides an isomorphism $\Br(S)\iso {\rm Tors}\,H^3(S,\Z)\cong \Z/2$. We will denote by $b_S$ the nonzero element of $\Br(S)$. 
\ind Let $\pi :X\rightarrow S$ be the $2$\tx to-1 cover of $S$ by a K3 surface. The aim of this note is to study the pull-back $\pi ^*b_S$ in $\Br(X)$.

\th Proposition
\enonce
 The class $\pi ^*b_S$ is represented, through the isomorphism \break $\Br(X)\iso\Hom(T_X,\Q/\Z)$, by the linear form $\tau \mapsto (\beta \cdot \pi _*\bar\tau )$, where $\bar\tau $ is the image of $\tau $ in $H^2(X,\Z/2)$ and $\beta $ any element of $H^2(S,\Z/2)$ which does not come from $H^2(S,\Z)$.
\endth 
{\it Proof} : Let $\beta $ be an element of $H^2(S,\Z/2)$ which does not come from $H^2(S,\Z)$, so that $p(\beta)=b_S\ (2.a) $. The pull-back $\pi ^*b_S\in\Br(X)$ is represented by $\pi ^*\beta \in H^2(X,\Z/2)$ $\cong H^2(X,\Z)\otimes \Z/2$;  its image in $\Hom(T_X,\Z/2)$ is the linear form
$\tau \mapsto (\pi ^*\beta \cdot \bar\tau )$. Since $(\pi ^*\beta \cdot \bar\tau )=(\beta \cdot \pi _*\bar\tau )$, the Proposition follows.\cqfd

\ind Part $(i)$ of the following Proposition shows that the class $\pi ^*\beta \in H^2(X,\Z/2)$ which appears above is nonzero. This does {\it not} say that $\pi ^*b_S$ is nonzero, as $\pi ^*\beta $ could come from a class in $\Pic(X)$ -- see \S 6.
 
\th Proposition
\enonce $(i)$ The kernel of $\pi ^*:H^2(S,\Z/2)\rightarrow H^2(X,\Z/2)$ is $\{0,k_S\}$.
\ind $(ii)$ The Gysin map $\pi _*:H^2(X,\Z)\rightarrow H^2(S,\Z)$ is surjective.
\endth
{\it Proof} :  To prove $(i)$ we
use the Hochschild-Serre spectral sequence :
$$E^{p,q}_2=H^p(\Z/2 , H^q(X,\Z/2)\ \Rightarrow \ H^{p+q}(S,\Z/2)\ .$$
 We have $E^{1,1}_2=0$, and $E^{2,0}_{\infty}=E^{2,0}_2=H^2(\Z/2 ,\Z/2)=\Z/2$. Thus the kernel of  $\pi ^*: $ $H^2(S,\Z/2)\rightarrow  H^2(X,\Z/2)$ is isomorphic to $\Z/2$. Since it contains $k_S$, it is equal to
$\{0,k_S\}$.
  \smallskip

  \ind Let us prove $(ii)$. Because of the formula $\pi _*\pi ^*\alpha =2\alpha $, the cokernel of $\pi _*:$ $H^2(X,\Z)\rightarrow H^2(S,\Z)$ is a $(\Z/2)$\tx vector space; therefore it suffices to prove that the transpose map
 $${}^t\pi _*:\Hom(H^2(S,\Z),\Z/2)\longrightarrow \Hom(H^2(X,\Z),\Z/2)$$ is injective.
 This is implied by the commutative diagram
 $$\xymatrix@M=4pt{H^2(S,\Z/2)\ar@{->>}[r]_<<<<<{v_S}\ar[d]^{\pi ^*}&\Hom(H^2(S,\Z),\Z/2)\ar[d]_{{}^t\pi _*}\\
                       H^2(X,\Z/2)\ar[r]^<<<<<{\sim}_<<<<<{v_X}&\Hom(H^2(X,\Z),\Z/2)
 }$$plus the fact that $\Ker\pi ^*=\Ker v_S=\{0,k_S\}$ (by $(i)$ and $(3.a)$).\cqfd
\section{Brauer groups and cyclic coverings}
\th Proposition
\enonce Let $\pi :X\rightarrow S$ be an \'etale, cyclic covering of smooth projective varieties over an algebraically closed field $k$. Let $\sigma $ be a generator of the Galois group $G$ of $\pi $, and let $\Nm:\Pic(X)\rightarrow \Pic(S)$ be the norm homomorphism. The kernel of  $\pi ^*:\Br(S)\rightarrow  \Br(X)$ is canonically isomorphic to $\Ker \Nm/(1-\sigma ^*)(\Pic(X))$. 
\endth
{\it Proof} : We
 consider the Hochschild-Serre spectral sequence $$E^{p,q}_2=H^p(G , H^q(X,\gm)\ \Rightarrow \ H^{p+q}(S,\gm)\ .$$
 
 Since $E_2^{2,0}=H^2(G ,k^*)=0$, the kernel of $\pi ^*:\Br(S)\rightarrow \Br(X)$ is identified with $E^{1,1}_{\infty}=\Ker \bigl(d_2: E^{1,1}_2\rightarrow E^{3,0}_2
\bigr)$. We have $E^{3,0}_2=H^3(G ,k^*)$; by periodicity of the cohomology of $G$, this group is canonically isomorphic to $H^1(G,k^*)=\Hom(G,k^*)$, the character group of $G$, which we denote by $\widehat{G}$.   So we view $d_2$ as a map from 
$H^1(G ,\Pic(X))$ to $\widehat{ G}$. 
\ind Let $\sym$ be the endomorphism $L\mapsto\bot_{g\in G}g^*L$ of $\Pic(X)$; recall that $H^1(G ,\Pic(X))$ is isomorphic to
$\Ker \sym/\im(1-\sigma^* )$.  We have $\pi ^*\Nm(L)=\sym(L)$ for $L\in\Pic(X)$, hence $\Nm$ maps $\Ker \sym$ into $\Ker \pi ^*\subset \Pic(S)$. 
 Now recall that $\Ker \pi ^*$ is canonically isomorphic to  $\widehat{ G}$: to  $\chi \in\widehat{ G}$ corresponds the subsheaf $L_\chi $ of $\pi _*{\cal O}_X$ where $G$ acts through the character $\chi $. Since $\Nm\rond(1-\sigma^* )=0$, the norm induces a homomorphism  $H^1(G ,\Pic(X))\rightarrow \Ker \pi ^*\cong  \widehat{G}$. The Proposition will follow from:
\th Lemma
\enonce The map $d_2 : H^1(G ,\Pic(X))\rightarrow \widehat{G}$ coincides with the homomorphism induced by the norm.
\endth
{\it Proof of the lemma} : We apply the formalism of [S], Proposition 1.1, where a very close situation is considered. This Proposition, together with property (1) which follows it, tells us that $d_2$ is given by cup-product with the extension class in $\Ext^2_G(\Pic(X),k^*)$  of the exact sequence of $G$\tx modules
$$1\rightarrow k^* \longrightarrow R_X^*\rightarrow \Div(X)\rightarrow \Pic(X)\rightarrow 0\ ,$$
where $R_X$ is the field of rational functions on $X$. This means that $d_2$ is the composition 
$$ H^1(G ,\Pic(X))\qfl{\partial}H^2(G, R_X^*/k^*)\qfl{\partial'}H^3(G,k^*)$$where $\partial$ and $\partial'$ are the 
coboundary maps associated to the short exact sequences
$$\nospacedmath\displaylines{0\rightarrow R_X^*/k^*\rightarrow \Div(X)\rightarrow \Pic(X)\rightarrow 0 \cr
\hbox{and}\hfill0\rightarrow k^*\rightarrow R_X^*\rightarrow R_X^*/k^*\rightarrow 0\ .\hfill
}$$

\ind Let $\lambda \in H^1(G ,\Pic(X))$, represented by $L\in\Pic(X)$ with $\bot _{g\in G}g^*L\cong {\cal O}_X$. Let $D\in\Div(X)$ such that $L={\cal O}_X(D)$. Then  $\sum_g g^*D$ is the divisor of a rational function $\psi \in R_X^*$, whose class  in $R_X^*/k^*$ is well-defined. This class is invariant under $G$, and defines the element  $\partial(\lambda )\in H^2(G ,R_X^*/k^*)$.
Since $\div \psi $ is invariant under $G$, there exists a character $\chi \in\widehat{G}$ such that $g^*\psi =\chi (g)\psi $ for each $g\in G$.
 Then $d^{1,1}_2(\lambda )=\chi  $ viewed as an element of $H^3(G,k^*)=\widehat{G}$. 
\ind It remains to prove that 
${\cal O}_S(\pi _*D ) =L_\chi $. Since $\div(\psi )=\pi ^*\pi_*D$, multiplication by $\psi $ induces a global isomorphism $u:\pi ^*{\cal O}_S(\pi _*D)\iso {\cal O}_X  $. Let $\varphi \in R_X$ be a  generator  of ${\cal O}_X(D) $ on an open $G$\tx invariant subset $U$ of $X$. Then $\Nm(\varphi )$ is a generator of ${\cal O}_S(\pi _*D)$ on $\pi (U)$, and $\pi ^*\Nm(\varphi ) $ is a generator of $\pi ^*{\cal O}_S(\pi _*D)$ on $U$; the function $h:=\psi\, \pi ^*\!\Nm(\varphi )$ on $U$ satisfies $g^*h=\chi (g)h$ for all $g\in G$. This proves that the homomorphism
$u^\flat:{\cal O}_S(\pi _*D)\rightarrow \pi _*{\cal O}_X $ deduced from $u$ maps  ${\cal O}_S(\pi _*D)$  onto the subsheaf $L_\chi $ of  $\pi _*{\cal O}_X$, hence our assertion.\cqfd

\ind We will need a complement of the Proposition in the complex case:
\th Corollary
\enonce Assume $k=\C$, and $H^1(X,{\cal O}_X )=H^2(S,{\cal O}_S )=0$.
The following conditions are equivalent:
\ind $(i)$ The map $\pi ^*:\Br(S)\rightarrow  \Br(X)$ is injective;
 \ind $(ii)$
Every class $\lambda= c_1(L)\in H^2(X,\Z)$, with  $L\in\Pic(X)$ and $\pi _*\lambda =0$, belongs to 
 $ (1-\sigma ^*) (H^2(X,\Z))$.   
\endth
\ind Observe that the hypotheses of the Corollary are satisfied when $S$ is a complex Enriques surfaces and $\pi :X\rightarrow S$ its universal cover.
\smallskip

{\it Proof} : By Proposition 4.1 $(i)$ is equivalent to $[L] =0$ in $H^1(G ,\Pic(X))$ for every $L\in\Pic(X)$ with $\Nm(L)= {\cal O}_S $, while $(ii)$ means $c_1(L)=0$ in $H^1(G ,H^2(X,\Z))$ for every such $L$. Therefore it suffices to prove that 
 the map 
$$H^1(c_1):H^1(G ,\Pic(X))\rightarrow  H^1(G ,H^2(X,\Z))$$ is injective.
\ind Since $H^1(X,{\cal O}_X )=0$ we have an exact sequence
$$0\rightarrow \Pic(X)\qfl{c_1}H^2(X,\Z)\longrightarrow Q\rightarrow 0\qquad \hbox{with }Q\subset H^2(X,{\cal O}_X )\ . $$Since $H^2(S,{\cal O}_S )=0$, there is no nonzero invariant vector in $H^2(X,{\cal O}_X )$, hence in $Q$. Then the associated long exact sequence implies that $H^1(c_1)$ is injective.\cqfd

\section{More topology}
\subsection As in \S 3, we denote by $S$ a complex Enriques surface and by $\pi :X\rightarrow S$ its universal cover. Thus $X$ is a K3 surface, with a fixed-point free involution $\sigma $ such that $\pi \rond \sigma =\pi $.
We will need some more precise results on the topology of the surfaces $X$ and $S$. 
We refer again to  [BHPV], ch. VIII.
\ind Let $E$ be the lattice $(-E_8)\oplus H$, where $H$ is the rank 2 hyperbolic lattice. 
 Let $H^2(S,\Z)_{\rm tf}$ be the quotient of $H^2(S,\Z)$ by its torsion subgroup $\{0,K_S\}$. We have isomorphisms
 $$H^2(S,\Z)_{\rm tf}\cong E \qquad H^2(X,\Z)\cong E\oplus E \oplus H$$
such that $\pi ^*:H^2(S,\Z)_{\rm tf}\rightarrow H^2(X,\Z)$ is identified with the diagonal embedding  $\delta : $ $E\hookrightarrow E\oplus E$, and $\sigma ^*$ is identified with the involution
$$\rho :(\alpha ,\alpha ',\beta )\mapsto (\alpha ',\alpha ,-\beta )\quad \hbox{of }\ E\oplus E \oplus H\ .$$\vskip-5pt
\subsection We consider now the cohomology with values in $\Z/2$. For a lattice $M$, we will write $M_2:=M/2M$. The scalar product of $M$ induces a product $M_2\otimes M_2\rightarrow \Z/2$; if moreover $M$ is {\it even}, there is a natural quadratic form $q:M_2\rightarrow \Z/2$ associated with that product, defined by $q(m)={1\over 2}\tilde m^2$, where $\tilde m\in M$ is any lift of $m\in M_2$.
In particular, $H_2$ contains a unique  element $\varepsilon $ with $q(\varepsilon )=1$: it is the class of $e+f$ where $(e,f)$ is a hyperbolic basis of $H$.

\ind  Using the previous isomorphism we identify  $H^2(X,\Z/2)$ with $E_2\oplus E_2 \oplus H_2$.
\th Proposition
\enonce The image of $\pi ^*:H^2(S,\Z/2)\rightarrow H^2(X,\Z/2)$ is $\delta(E_2)\oplus
(\Z/2)\varepsilon $.
\endth
{\it Proof} : This image is invariant under $\sigma^* $, hence is contained in $\delta (E_2) \oplus H_2$; by Proposition $3.6 \ (i)$ it is $11$\tx dimensional, hence a hyperplane in $\delta (E_2)  \oplus H_2$,  containing $\delta (E_2)$ (which is spanned by the classes coming from $H^2(S,\Z)$). So $\pi ^*H^2(S,\Z/2)$ is spanned by $\delta (E_2)$ and a
 nonzero element of $H_2$;  it suffices to prove that this element is $\varepsilon $. Since the elements  of $H^2(S,\Z/2)$ which do not come from $H^2(S,\Z)$ have square 1 (3.2), this is a consequence of the following lemma.\cqfd 

\th Lemma
\enonce For every $\alpha \in H^2(S,\Z/2)$, $q(\pi ^*\alpha )=\alpha ^2$.
\endth
{\it Proof} : This  proof has been shown to me by J. Lannes. The key ingredient is the {\it Pontryagin square}, a cohomological operation
$${\cal P}:H^{2m}(M,\Z/2)\longrightarrow H^{4m}(M,\Z/4)$$
defined for any reasonable topological space $M$ and satisfying
 a number of interesting properties (see [M-T], ch. 2, exerc. 1). We will state only those we  need in the case of interest for us, namely $m=2$ and $M$ is a compact oriented 4-manifold. We identify $H^4(M,\Z/4)$ with $\Z/4$; then ${\cal P}:H^2(M,\Z)\rightarrow \Z/4$ satisfies:
\ind a) For $\alpha \in H^2(M,\Z/2)$, the class of ${\cal P}(\alpha )$ in $\Z/2$ is $\alpha ^2$;
\ind b) If $\alpha  \in H^2(M,\Z/2)$ comes from $\tilde\alpha  \in H^2(M,\Z)$, then
${\cal P}(\alpha )=\tilde\alpha^2 $ (mod. 4). 
 In particular, if $M$ is a K3 surface, we have ${\cal P}(\alpha )=2q(\alpha )$ in $\Z/4$. 

\ind Coming back to our situation, let $\alpha \in H^2(S,\Z/2)$. We have in $\Z/4$:
$$\eqalign{{\cal P}(\pi ^*\alpha )&= 2\,{\cal P}(\alpha )\quad&\hbox{by functoriality}\cr
&=2\,\alpha ^2\quad&\hbox{by a), and}\cr
{\cal P}(\pi ^*\alpha )&=2\,q(\pi ^*\alpha )\quad&\hbox{by b).}
}$$Comparing the two last lines gives  the lemma.\cqfd
\th Corollary
\enonce The kernel of $\pi _*:H_2\rightarrow \{0,k_S\}$ is $\{0,\varepsilon \}$.
\endth
{\it Proof} : By Proposition 5.3 $\varepsilon $ belongs to $\im \pi ^*$, hence $\pi _*\varepsilon =0$. It remains to check that $\pi _*$ is nonzero on $H^1(\Z/2 ,H^2(X,\Z))\cong H_2$. We know that there is an element $\alpha \in H^2(X,\Z)$ with $\pi _*\alpha =K_S$ (Prop. 3.6 $(ii)$); it belongs to $\Ker(1+\sigma ^*)$, hence defines an element $\bar\alpha $ of $H^1(\Z/2 ,H^2(X,\Z))$ with $\pi _*\bar\alpha \neq 0$.\cqfd

\th Corollary
\enonce Let $\lambda \in H^2(X,\Z)$. The following conditions are equivalent:
\ind $(i)$  $\pi _*\lambda =0$ and $\lambda \notin (1-\sigma ^*) (H^2(X,\Z))\,;$
\ind $(ii)$ $\sigma ^*\lambda =-\lambda $ and $\lambda ^2\equiv 2$ {\rm (mod. 4)}. 
\endth
{\it Proof} : Write $\lambda =(\alpha ,\alpha ',\beta )\in E\oplus E\oplus H$; let $\bar\beta $ be the class of $\beta $ in $H_2$. Both conditions imply $\sigma ^*\lambda =-\lambda $, hence $\alpha '=-\alpha $. Since $(\alpha ,-\alpha )=(1-\sigma ^*)(\alpha ,0)$ and $2\beta =(1-\sigma ^*)(\beta )$, the conditions of $(i)$ are equivalent to  $\pi _*\bar\beta =0$ and $\bar\beta  \neq 0$, that is, $\bar\beta =\varepsilon $ (Corollary 5.5). On the other hand we have $\lambda ^2=2\alpha ^2+\beta ^2\equiv 2q(\bar\beta )$ (mod.~4), hence $(ii)$ is also equivalent to $\bar\beta =\varepsilon $.\cqfd

\ind We can thus rephrase Corollary 4.3 in our situation:

\th Corollary
\enonce We have $\pi ^*b_S=0$ if and only if there exists a line bundle $L$ on $X$ with $\sigma ^*L=L^{-1}$ and $c_1(L)^2\equiv 2$ {\rm (mod. 4)}.\cqfd
\endth
\bigskip

{\eightpoint\baselineskip=12pt
\leftskip1cm\rightskip1cm
\rem{Remark} My original proof of (5.3-5) was less direct and less general, but still  perhaps of some interest. The key point is to show that on $H_2$ $q$ takes the value 1 exactly on the  nonzero element of $\Ker \pi _*$, or equivalently that an element $\alpha \in H_2$ with $\pi^{} _*\alpha =k_S$ satisfies $q(\alpha )=0$. Using deformation theory (see (6.1) below), one can assume that $\alpha $ comes from a class in $\Pic(X)$. To conclude I applied the following lemma:
\th Lemma
\enonce  Let $L$ be a line bundle on $X$ with $\Nm(L)=K_S $. Then $c_1(L)^2$ is divisible by $4$.
\endth
{\it Proof} : Consider the rank 2 vector bundle $E=\pi _*(L)$. The norm induces a non-degenerate quadratic form $N:{\rm Sym}^2E\rightarrow K_S$
([EGA2], 6.5.5). In particular, $N$ induces an isomorphism
   $E\iso E^*\otimes K_S$, and defines a pairing
   $$H^1(S,E)\otimes H^1(S,E)\rightarrow H^2(S,K_S)\cong {\bb C}$$
which is alternating and non-degenerate. Thus $h^1(E)$ is even; since $h^0(E)=h^2(E)$ by Serre duality, $\chi (E)$ is even, and so is $\chi (L)=\chi (E)$. By Riemann-Roch this implies that ${1\over 2}c_1(L)^2$ is even.\cqfd}

\font\gragrec=cmmib10 
\def\pig{\hbox{\gragrec \char25}}
\def\bg{\hbox{\gragrec \char98}}

\section{The vanishing of ${\bf \pig ^*\bg_S}$ on the moduli space}
\subsection 
We briefly recall the theory of the period map for Enriques surfaces, due to Horikawa (see [BHPV], ch. VIII, or [N]). We keep the notations of (5.1). We denote by $L$ the lattice $E\oplus E\oplus H$, and by $L^-$  the $(-1)$\tx eigenspace of  the involution $\rho:(\alpha ,\alpha ',\beta )\mapsto (\alpha ',\alpha ,-\beta )$, that is, the submodule of elements $(\alpha ,-\alpha ,\beta )$.
\ind A {\it marking} of the Enriques surface $S$ is an isometry
$\varphi : H^2(X,\Z)\rightarrow L$ which conjugates $\sigma ^*$ to  $\rho $. 
The line  $H^{2,0}\subset H^2(X,\C)$ is anti-invariant under $\sigma ^*$, so its image by
$\varphi _{\C}:H^2(X,\C)\rightarrow L_{\C}$ lies in $L_{\C}^-$. The corresponding point $[\omega ]$ of $\P(L_{\C}^-)$ is  the {\it period}  $\wp(S,\varphi )$. It belongs to the domain $\Omega \subset \P(L_{\C}^-)$ defined by the equations
$$(\omega \cdot \omega )=0\qquad (\omega \cdot \bar\omega )>0\qquad (\omega \cdot \lambda )\neq 0\quad \hbox{for all }\lambda \in L^-\hbox{ with }\lambda ^2=-2\ .$$
This is an analytic manifold, which is the moduli space for marked Enriques surfaces.
To each class $\lambda \in L^-$ we associate the hypersurface $H_\lambda $ of $\Omega $ defined by $(\lambda \cdot \omega )=0$.

\th Proposition
\enonce We have $\pi ^*b_S=0$ if and only if $\wp(S,\varphi )$ belongs to one of the hypersurfaces $H_\lambda $ for some vector $\lambda \in L^-$ with $\lambda ^2\equiv 2\ {\rm (mod.\ 4)} $.
\endth
{\it Proof} :  The period point $\wp(S,\varphi )$ belongs to $H_\lambda $ if and only if  $\lambda $ belongs to $c_1(\Pic(X))$; by Corollary 5.7, this is equivalent to $\pi ^*b_S=0$.\cqfd

\ind To get a complete picture we want to know which of the $H_\lambda $ are really needed:

\th Lemma 
\enonce Let $\lambda $ be a primitive element of $L^-$.
\ind $(i)$ The hypersurface $H_\lambda $ is non-empty if and only if $\lambda ^2<-2$.
\ind $(ii)$ If $\mu $ is another primitive element of $L^-$ with 
$H_\mu =H_\lambda\neq \varnothing$, then $\mu =\pm \lambda $.
\endth
{\it Proof} : Let $W $ be the subset of $L^-_{\C}$ defined by the conditions $\omega ^2=0$, $\omega \cdot \bar\omega >0$. If
we write $\omega =\alpha +i\beta $ with $\alpha ,\beta \in L^-_{\R}$, these conditions translate as $\alpha ^2=\beta ^2>0$, $\alpha \cdot \beta =0$. Thus  $  W  \cap \lambda ^\perp\neq \varnothing$ 
is equivalent to the existence of a positive 2-plane in $L^-_{\R}$ orthogonal to $\lambda $. Since $L^-$ has signature $(2,10)$, this is  also equivalent to $\lambda ^2<0$. 
\ind If $W \cap \lambda ^\perp$ is non-empty, $\lambda ^\perp $ is the only hyperplane containing it, and  $\C\,\lambda $ is the orthogonal of $\lambda ^\perp $ in $L^-$. Then $\lambda $ and $-\lambda $ are the only primitive vectors of $L^-$ contained in $\C\,\lambda$. In particular  $\lambda $ is determined up to sign by $H_\lambda $, which proves $(ii)$.
\ind  Let us prove $(i)$. We have seen that $H_\lambda $ is empty for $\lambda ^2\geq 0$, and also for $\lambda ^2=-2$ by definition of $\Omega $. Assume $\lambda ^2<-2$ and $H_\lambda =\varnothing$; then $H_\lambda $ must be contained in one of the hyperplanes $H_\mu $ with $\mu ^2=-2$; by $(ii)$  this implies $\lambda =\pm \mu $, a contradiction.\cqfd

\subsection   Let $\Gamma $ be the group of isometries of $L^-$.  The group $\Gamma $ acts properly discontinuously on $\Omega $, and the quotient ${\cal M}=\Omega /\Gamma $ is a quasi-projective variety. 
The image  in ${\cal M}$ of the period $\wp(S,\varphi )$ does not depend on the choice of $\varphi $; let us denote it by $\wp(S)$.  The map $S\mapsto \wp(S)$ induces a bijection between isomorphism classes of Enriques surfaces and ${\cal M}$; the variety ${\cal M}$ is a
 (coarse) moduli space for Enriques surfaces. 
\th Corollary
\enonce The surfaces $S$ for which $\pi ^*b_S=0$ form an infinite, countable  union of (non-empty) hypersurfaces in the moduli space ${\cal M}$.
\endth
{\it Proof} : Let $\Lambda $ be the set of primitive elements $\lambda$ in $L^-$ with $\lambda ^2<-2$ and $\lambda ^2\equiv 2\ {\rm (mod.\ 4)} $. For $\lambda \in\Lambda $, let ${\cal H}_\lambda $ be the image of $H_\lambda $ in ${\cal M}$;  the argument of [BHPV], ch. VIII, Cor. 20.7 shows that ${\cal H}_\lambda $ is an algebraic hypersurface in ${\cal M}$. 
By Proposition 6.2 and Lemma 6.3 the surfaces $S$ with $\pi ^*(b_S)=0$ form the subset $\bigcup_{\lambda\in\Lambda } {\cal H}_\lambda$. By Lemma 6.3 $(ii)$ we have ${\cal H}_\lambda ={\cal H}_\mu $ if and only if $\mu =\pm g\lambda $ for some element $g$ of $\Gamma $. This implies $\lambda ^2=\mu ^2$; but $\lambda ^2$ can 
be any  number of the form $-2k$ with $k$ odd $>1$ (take for instance $\lambda =e-kf$, where $(e,f)$ is a hyperbolic basis of $H$), so there are infinitely many distinct hypersurfaces among the ${\cal H}_\lambda $.\cqfd

\vskip2truecm
\font\cc=cmcsc10

\centerline{ REFERENCES} \vglue15pt\baselineskip12.8pt
\def\num#1#2#3{\smallskip\item{\hbox to\parindent{\enskip [#1]\hfill}}{\cc #2: }{\sl #3}}
\parindent=1.3cm 
\num{BHPV}{\ W. Barth,  K. Hulek, C. Peters, A. Van de Ven} {Compact complex surfaces}. 2nd edition. Ergebnisse der Mathematik und ihrer Grenzgebiete, {\bf 4}. Springer-Verlag, Berlin (2004).
\num{EGA2}{\ A.~Grothendieck}{\'El\'ements de g\'eom\'etrie alg\'ebrique} II.  Publ. Math. IHES {\bf 8} (1961).
\num{G}{A. Grothendieck}{Le groupe de Brauer} I--II.  Dix Expos\'es sur la Cohomologie des Sch\'emas,  pp. 46--87; North-Holland, Amsterdam (1968).
\num{H-S}{D. Harari, A. Skorobogatov} {Non-abelian descent and the arithmetic of Enriques surfaces}. Intern. Math. Res. Notices  {\bf 52} (2005), 3203--3228.

\num{M-S}{J. Milnor, J. Stasheff}{Characteristic classes}. Annals of Math. Studies {\bf 76}. Princeton University Press, Princeton  (1974).
\num{M-T}{R. Mosher, M. Tangora}{Cohomology operations and applications in homotopy theory}. Harper \& Row, New York-London (1968).
\num{N}{Y. Namikawa}{Periods of Enriques surfaces}.  Math. Ann.  {\bf 270}  (1985),  no. 2, 201--222.
\num{S}{A. Skorobogatov}
{On the elementary obstruction to the existence of rational points}. 
Math. Notes {\bf 81} (2007), no. 1, 97--107.
\vskip1cm
\def\pc#1{\eightrm#1\sixrm}
\hfill\vtop{\eightrm\hbox to 5cm{\hfill Arnaud {\pc BEAUVILLE}\hfill}
  \hbox to 5cm{\hfill Laboratoire J.-A. Dieudonn\'e\hfill}
 \hbox to 5cm{\sixrm\hfill UMR 6621 du CNRS\hfill}
\hbox to 5cm{\hfill {\pc UNIVERSIT\'E DE}  {\pc NICE}\hfill}
\hbox to 5cm{\hfill  Parc Valrose\hfill}
\hbox to 5cm{\hfill F-06108 {\pc NICE} Cedex 02\hfill}}
\end